\documentclass[oneside,english,reqno]{amsart}
\usepackage{palatino}
\usepackage[T1]{fontenc}
\usepackage[latin1]{inputenc}
\usepackage{float}
\usepackage{amssymb}
\IfFileExists{url.sty}{\usepackage{url}}
                      {\newcommand{\url}{\texttt}}

\makeatletter

\providecommand{\tabularnewline}{\\}

  \theoremstyle{plain}
  \newtheorem*{thm*}{Theorem}
  \theoremstyle{plain}
  \newtheorem{lem}{Lemma}
  \theoremstyle{remark}
  \newtheorem*{rem*}{Remark}
 \theoremstyle{definition}
 \newtheorem*{defn*}{Definition}

\DeclareMathOperator{\spec}{spec}

\DeclareMathOperator{\wind}{wind}
\DeclareMathOperator{\Osc}{Osc}
\def\esssup{\mathop{\rm ess\,sup}}
\newcommand{\VMO}{\mathrm{VMO}}
\newcommand{\BMO}{\mathrm{BMO}}

\ifx\hyperlink\undefined

\newcommand{\hyperlink}[2]{#2}
\else

\fi
\usepackage{graphicx}
\usepackage{color}

\makeatother
\begin{document}

\title{One cannot hear the winding number}

\author{Jean Bourgain}

\address{JB: Institute for Advanced Study, 1 Einstein dr., Princeton NJ 08540
USA}

\email{bourgain@ias.edu}

\author{Gady Kozma}

\thanks{This material is based upon work supported by the National Science
Foundation under agreement DMS-0111298. Any opinions, findings and
conclusions or recommendations expressed in this material are those
of the authors and do not necessarily reflect the views of the National
Science Foundations}

\address{GK: The Weizmann Institute of Science, Rehovot 76100, Israel.}

\email{gady.kozma@weizmann.ac.il}

\begin{abstract}
We construct an example of two continuous maps $f$ and $g$ of the
circle to itself with $|\widehat{f}(n)|=|\widehat{g}(n)|$ for all
$n\in\mathbb{Z}$ but with different winding numbers, answering a
question of Brezis.
\end{abstract}
\maketitle

\section{Introduction}

A continuous cycle $f$ in $\mathbb{C}\setminus\{0\}$ has a well
defined winding number around zero, which we shall denote by $\wind f$.
If $f$ is smooth and its image is in the circle $\mathbb{T}=\{ z:|z|=1\}$
then the winding number has an elegant formula using the Fourier coefficients
of $f$ (we consider $f$ as a function from $[0,1]$ to $\mathbb{T}$).
Indeed, by the residue formula for the function $z^{-1}$ we have\[
\wind f=\frac{1}{2\pi i}\int\frac{f'}{f}\]
and since $1/f=\overline{f}$ we get\begin{equation}
\wind f=\frac{1}{2\pi i}\int f'\overline{f}=\frac{1}{2\pi i}\sum\widehat{f'}(n)\overline{\widehat{f}(n)}=\sum n|\widehat{f}(n)|^{2}.\label{eq:wind Fourier}\end{equation}
This paper is part of a line of research trying to understand what
role does smoothness play in (\ref{eq:wind Fourier}). See Brezis
\cite{B} for a fascinating review of related results and problems,
high dimensional analogs and applications.

The earliest investigations in this direction are due to L. Boutet
de Monvel and O. Gabber. They realized that the left-most equality
in (\ref{eq:wind Fourier}) makes sense for functions in the fractional
Sobolev space $W^{1/2,2}$, when the integral is understood in the
sense of $W^{1/2,2}$--$W^{-1/2,2}$ duality (see the appendix of
\cite{BGP91}). This equality allowed them to extend the notion of
winding number to the discontinuous part of $W^{1/2,2}$. By approximating
with continuous functions they showed that this generalized winding
number was still an integer. Interestingly, $W^{1/2,2}$ is exactly
the space of functions for which the right hand side of (\ref{eq:wind Fourier})
converges absolutely, but apparently Boutet de Monvel and Gabber were
not aware of (\ref{eq:wind Fourier}). The connection to Fourier expansion
was discovered by Brezis in 1995 (following a question of Gelfand)
and immediately begot new questions.

As a side remark, defining the winding number using approximation
by continuous functions is most natural in the space $\VMO$ of functions
of vanishing mean oscillation ($\VMO$ is the closure of the continuous
functions in $\BMO$ and, heuristically, it relates to $\BMO$ like
continuous functions relate to $L^{\infty}$). Note that $\VMO$ contains
$W^{1/2,2}$ and the $\VMO$-winding number agrees with both the definition
of \cite{BGP91} and with (\ref{eq:wind Fourier}) on $W^{1/2,2}$,
see \cite{B}. A high dimensional analog, the $\VMO$-degree, was
developed by Brezis and Nirenberg \cite{BN95}, see also \cite{BC83}.

Returning to (\ref{eq:wind Fourier}), when one leaves $W^{1/2,2}$
the picture gets more complicated. Following a question of Brezis,
Korevaar \cite{K99} showed that in general the sum (\ref{eq:wind Fourier}),
considered as the balanced limit\[
\lim_{N\to\infty}\sum_{n=-N}^{N}n|\widehat{f}(n)|^{2}\]
 may converge to any desired value different from the winding number,
including $\pm\infty$, or may diverge. Replacing convergence with
Abel summability does not change the picture. On the other hand Korevaar
shows that (\ref{eq:wind Fourier}) does hold for continuous functions
with bounded variation (this class is not contained in $W^{1/2,2}$).
More intriguing, perhaps, is Kahane's proof \cite{K05} that for H\"older
functions with exponent $>\frac{1}{3}$ one can retrieve the winding
number by summing (\ref{eq:wind Fourier}) using Riemann's summation
method. Brezis \cite[theorem 5]{B} then noticed that Kahane's proof
works in the space $W^{1/3,3}$. In particular Kahane's argument shows
that in $W^{1/3,3}$ the absolute values of the Fourier transform
determine the winding number, which can be picturesquely described
as {}``hearing the winding number''. It is conjectured that one
can hear the winding number in the class $W^{1/p,p}$ for any $p$.
Here $W^{1/p,p}$ is in the sense of Sobolev-Slobodecki\u\i{} spaces,
see e.g.~\cite{KJF77}. We remind the reader that in this case the
definition of these spaces simplifies to the condition \[
W^{1/p,p}=\left\{ f:\int\int\frac{|f(x)-f(y)|^{p}}{|x-y|^{2}}<\infty\right\} .\]
This is an increasing family of classes  and all $W^{1/p,p}\subset\VMO$. 

Korevaar's negative result, while strongly hinting that one cannot
hear the winding number in general does not actually preclude it.
In this paper we show that some kind of smoothness must be assumed
in order to hear the winding number. Specifically we show

\begin{thm*}
There exists two continuous functions $f,g:\mathbb{T}\to\mathbb{T}$
with $|\widehat{f}(n)|=|\widehat{g}(n)|$ for all $n$ and different
winding numbers.
\end{thm*}
The constructed $f$ and $g$ are highly singular. We made no attempt
to optimize the construction in this respect nor do we believe that
our techniques are adequate for this purpose. In fact, even forcing
$f$ and $g$ to be continuous, as opposed to just $\VMO$, required
a non-negligible increase in the complexity of the proof. Another
interesting question which we cannot address at this is time is whether
one can strengthen $|\widehat{f}(n)|=|\widehat{g}(n)|$ to $\widehat{f}(n)=\pm\widehat{g}(n)$.

We remark that functions with identical absolute value and identical
absolute value of the Fourier transform are called \emph{Pauli partners}.
Thus in this language the theorem states that there exists two Pauli
partners $f$ and $g$ with $|f|=1$ and different winding numbers.
For many interesting techniques for producing Pauli partners, see
\cite{J99}.

\subsection{About the proof}

Inspired by de Leeuw-Kahane-Katznelson \cite{dLKK77} we construct
$f$ and $g$ by an iterative correction scheme with each stage combining
a deterministic step and a probabilistic step. We start with functions
which take values outside $\mathbb{T}$, so (\ref{eq:wind Fourier})
does not hold and one may construct trigonometric polynomials with
$|\widehat{f}(n)|=|\widehat{g}(n)|$ and different winding numbers
quite easily (see page \pageref{page:start}). At each stage the polynomials
are extended so as to continue satisfying $|\widehat{f}(n)|=|\widehat{g}(n)|$
and preserve their winding number while getting closer and closer
to $\mathbb{T}$.

Let us describe one correction. Assume $|f|$ is too low in an interval
$I$. We always correct upwards, increasing $|f|$ --- this is a corollary
of the fact that extending a polynomial increases its $L^{2}$ norm
--- so our lemmas are usually formulated non-symmetrically with respect
to $|f|$. We wish to correct $f$ on $I$ (call the correction $F$)
such that $F$ has bigger absolute value on $I$ but is highly oscillatory
such that $F-f$ lives essentially in the high end of the spectrum.

We employ two different techniques to get $F$ from $f$. The first,
and simpler, is described in lemma \ref{lem:arcsin}. Its advantage
is that $F$ has the desired properties on all of $I$. Its disadvantage
is that $F-f$ has no structure that can be used in order to extend
the other polynomial, $g$. Hence when using this technique we simply
define \[
G=g+\sum\pm\widehat{F-f}(n)e^{int}\]
where $\pm$ are random signs. If $I$ is sufficiently small then
$\left\Vert F-f\right\Vert _{2}$ would be small, and we could bound
$\left\Vert G-g\right\Vert _{\infty}$ efficiently by known properties
of random trigonometric series

The second technique is used in lemma \ref{lem:smooth} and especially
in lemma \ref{lem:walshIJ}. We construct $F-f$ so as to have many
small pieces sitting in different areas of the spectrum and use this
structure in order to construct $G$ so that its distance from a unimodular
function is decreased as little as possible. The disadvantage of this
technique is that we haven't figured out how to correct on the whole
of $I$. $F$ has the desired absolute value on most of $I$ but a
small exceptional set remains, and is handled using the first technique.

\subsection{Notations}

We shall consider functions from $[0,1]$ into $\mathbb{C}$, which
we consider as periodic e.g.~when we say that $f$ is continuous
we mean also that $f(0)=f(1)$. A \textbf{rotation} is the action
of transforming $f(t)$ into $f(t+a)$ for some $a\in[0,1]$, again
periodically. The Fourier coefficients are defined by $\widehat{f}(n)=\int_{0}^{1}f(t)e^{-2\pi int}$.
For a trigonometric polynomial $f$ we define its spectrum, denoted
by $\spec f$ and its degree, denoted by $\deg f$ by\[
\spec f:=\{ n\in\mathbb{Z}:\widehat{f}(n)\neq0\}\qquad\deg f:=\max\{|n|:n\in\spec f\}.\]
If $f$ and $F$ are two trigonometric polynomials we say that $F$
\textbf{extends} $f$ if they are identical on the hull of the spectrum
of $f$ namely\[
\widehat{F}(n)=\widehat{f}(n)\quad\forall|n|\leq\deg f.\]
$\omega(\delta;f)$ will denote the modulus of continuity, namely\[
\omega(\delta;f):=\max_{|x-y|\leq\delta}|f(x)-f(y)|.\]

We shall denote by $C^{n}$, $n\in\mathbb{N}$ the space of functions
with $n$ continuous derivative and by $C^{\alpha}$, $\alpha\in\left]0,1\right[$
the space of H\"older continuous functions of order $\alpha$. The
respective norms are\begin{align*}
\left\Vert f\right\Vert _{C^{n}} & :=\max\Vert f\Vert_{\infty},\Vert f^{(n)}\Vert_{\infty} & \left\Vert f\right\Vert _{C^{\alpha}} & :=\max_{x\neq y}\frac{|f(x)-f(y)|}{|x-y|^{\alpha}}\end{align*}
(we shall not need other $C^{x}$-s --- in fact we shall only use
$C^{2}$, $C^{1}$ and $C^{1/2}$). Other norms we will use are \begin{align*}
\left\Vert f\right\Vert _{2} & :=\sqrt{\int_{0}^{1}|f|^{2}} & \left\Vert f\right\Vert _{\infty} & :=\esssup_{t\in[0,1]}|f(t)|\end{align*}
and we denote by $V(f)$ the total variation of $f$. We also remind
the reader that \begin{equation}
V(fg)\leq V(f)\left\Vert g\right\Vert _{\infty}+V(g)\left\Vert f\right\Vert _{\infty}\quad\left\Vert fg\right\Vert _{C^{\alpha}}\leq\left\Vert f\right\Vert _{C^{\alpha}}\left\Vert g\right\Vert _{\infty}+\left\Vert g\right\Vert _{C^{\alpha}}\left\Vert f\right\Vert _{\infty}.\label{eq:prod}\end{equation}

By $f_{+}$ we mean $\max\{ f,0\}$. For a set $E$ we denote by $\mathbf{1}_{E}$
the indicator function which is $1$ on $E$ and $0$ elsewhere. $\left\lfloor x\right\rfloor $
and $\left\lceil x\right\rceil $ denote, respectively the floor and
ceil functions i.e.~the largest integer $\leq x$ and the smallest
integer $\geq x$. By $C$ and $c$ we shall denote absolute constants
whose precise value is unimportant as far as this paper is concerned,
and could change from formula to formula or even within the same formula.
$C$ will pertain to constants which are {}``big enough'' and $c$
to constants which are {}``small enough''. We will number a few
$C$ and $c$-s --- only those which we will reference later on. When
we say {}``$x$ is sufficiently large'' we mean {}``$x>C$''.

\section{Proof}

\subsection{The local correction scheme}

\begin{lem}
\label{lem:C2}Let $h$ be a $C^{2}$ function. Then\[
V(\sqrt{h_{+}})\leq\sqrt{\left\Vert h\right\Vert _{C^{2}}}.\]

\end{lem}
\begin{proof}
In fact we will prove the stronger \[
V(\sqrt{|h|})\leq\sqrt{\left\Vert h\right\Vert _{C^{2}}}.\]
Let $I=[a,b]$ be an interval such that $h'(a)=h'(b)=0$ and $h$
is monotone on $[a,b]$. Then \[
V_{I}(h)=|h(b)-h(a)|\leq(b-a)\max_{t\in I}|h'(t)|\leq{\textstyle \frac{1}{2}}(b-a)^{2}\max_{t\in I}|h''(t)|\]
where $V_{I}(h)$ is the variation of $h$ on $I$. If in addition
$h(t)\neq0$ for all $t\in\left]a,b\right[$ then \[
V_{I}(\sqrt{|h|})=\left|\sqrt{|h(b)|}-\sqrt{|h(a)|}\right|\leq\sqrt{|h(b)-h(a)|}\leq\sqrt{{\textstyle \frac{1}{2}}\left\Vert h\right\Vert _{C^{2}}}(b-a).\]
If $h(t)=0$ for some $t\in\left]a,b\right[$ then \[
V_{I}(\sqrt{|h|})=\sqrt{|h(b)|}+\sqrt{|h(a)|}\leq\sqrt{2|h(b)-h(a)|}\]
so in both cases\begin{equation}
V_{I}(\sqrt{|h|})\leq\sqrt{\left\Vert h\right\Vert _{C^{2}}}(b-a).\label{eq:VI}\end{equation}
A similar argument shows that if $h'(t)=0$ for any $t\in[a,b]$ (and
without assuming monotonicity of $h$) then \begin{equation}
\left|\sqrt{|h(a)|}-\sqrt{|h(b)|}\right|\leq\sqrt{\left\Vert h\right\Vert _{C^{2}}}(b-a).\label{eq:Vab}\end{equation}

Let now $a_{1}<a_{2}<\dotsb<a_{N}\in[0,1]$. We need to estimate the
variation with respect to $a_{1},\dotsc,a_{N}$. Clearly we may add
points and we add, for any segment $[a_{i},a_{i+1}]$ where $h$ is
not monotonic the maximal and minimal points in $[a_{i},a_{i+1}]$
where $h'=0$, removing duplicates if they arise (each may be equal
to the boundary point, and the two points may be equal). Denote the
new list also by $a_{i}$. Let $i_{1}<\dotsb<i_{K}$ the points where
$h'(a_{i_{k}})=0$ and assume that $a_{i_{1}}=0$ and $a_{i_{K}}=1$
(as we may, by rotating $h$ and the $a_{i}$-s and adding one $a_{i}$,
if necessary). It is now easy to verify that the following holds\begin{equation}
h\mbox{ is monotone on }[a_{i_{k}},a_{i_{k+1}}]\mbox{ whenever }i_{k+1}>i_{k}+1.\label{eq:monor0}\end{equation}
Let us now write\[
V(\sqrt{|h|};a_{1},\dotsc,a_{N}):=\sum_{i=1}^{N-1}\left|\sqrt{|h(a_{i+1})|}-\sqrt{|h(a_{i})|}\right|=\sum_{k=1}^{K-1}v_{k}\]
where \[
v_{k}:=\sum_{i=i_{k}}^{i_{k+1}-1}\left|\sqrt{|h(a_{i+1})|}-\sqrt{|h(a_{i})|}\right|.\]
However, by (\ref{eq:monor0}) we can use (\ref{eq:VI}) for the case
that $h$ is monotone on $[a_{i_{k}},a_{i_{k+1}}]$ and (\ref{eq:Vab})
in the case that $i_{k+1}=i_{k}+1$ and in both cases we get $v_{k}\leq\sqrt{\left\Vert h\right\Vert _{C^{2}}}(a_{i_{k+1}}-a_{i_{k}})$.
This finishes the lemma.
\end{proof}
\begin{lem}
\label{lem:arcsin}Let $f$ be a trigonometric polynomial with $1-c_{1}<|f|<1+c_{1}$
for some absolute constant $0<c_{1}<1$. Let $[a,b]\subset\mathbb{T}$
be some interval such that $|f(a)|=|f(b)|$ and $|f(t)|<|f(a)|$ for
all $t\in[a,b]$. Let $\epsilon\in\left]0,1\right[$ be some parameter.
Then one can extend $f$ as $F$ such that\begin{alignat}{2}
 & |\,|F(t)|-|f(a)|\,|<\epsilon &  & \forall t\in[a,b]\label{eq:lem|F|-|f|}\\
 & |F(t)-f(t)|<\epsilon &  & \forall t\not\in[a,b]\label{eq:F-f}\\
 & \deg F\leq C(\deg f)^{12}/\epsilon^{4} & \quad\label{eq:degF}\end{alignat}

\end{lem}
\begin{proof}
Assume $a=0$ (otherwise we can rotate the whole thing). We shall
need the following function, defined on $\left[0,\infty\right[$:\[
\varphi(x)=\sum_{j=1}^{\infty}\frac{(-1)^{j+1}}{(2j)!}\binom{2j}{j}4^{-j}x^{j}.\]
Clearly this is an analytic function with $\varphi(0)=0$ and $\varphi'(0)>0$.
Hence we may invert it in some neighborhood of $0$. We get an analytic
monotone function $\psi:[0,c]\to[0,1]$ with $\psi'\leq C$ and $\psi''\leq C$.

The construction now goes as follows. Denote $\tau=|f(0)|$, $N=\deg f$
and let $M=M(N)$ be some number to be fixed later. Define\begin{align*}
\delta(t) & =\begin{cases}
\sqrt{\psi\left(1-\frac{|f(t)|}{\tau}\right)} & t\in[0,b]\\
0 & \textrm{otherwise}\end{cases}\\
f_{2}(t) & =\begin{cases}
\frac{f(t)}{|f(t)|}\tau e^{i\delta(t)\sin Mt} & t\in[0,b]\\
f(t) & \textrm{otherwise}\end{cases}.\end{align*}
If $c_{1}$ is sufficiently small then $\delta$ is always well defined.
Fix some value of $c_{1}<\frac{1}{2}$ satisfying this requirement.
The following properties of $f_{2}$ now follow:\begin{alignat}{2}
|f_{2}(t)| & =\tau & \quad & \forall t\in[0,b]\label{eq:F>}\\
f_{2}(t) & =f(t) & \quad & \forall t\not\in[0,b].\label{eq:|F|-|f|}\end{alignat}

To estimate $\widehat{f_{2}}$ develop $e^{i\delta(t)\sin Mt}$ in
a Taylor expansion. We get\[
e^{i\delta(t)\sin Mt}=1+\sum_{j=1}^{\infty}\frac{1}{j!}(i\delta(t)\sin Mt)^{j}.\]
Expanding $\sin^{j}Mt=\left((e^{iMt}-e^{-iMt})/2i\right)^{j}$ we
get for odd $j$ a sum with no constant coefficient and for even $j$
the constant coefficient is $\binom{j}{j/2}2^{-j}$. Hence we may
write, for all $t\in[0,b]$,\begin{align*}
e^{i\delta(t)\sin Mt} & =1+\sum_{j=1}^{\infty}\frac{(-1)^{j}}{(2j)!}\binom{2j}{j}4^{-j}\delta^{2j}(t)+R(t)=1-\varphi(\delta^{2}(t))+R(t)=\\
 & =\frac{|f(t)|}{\tau}+R(t)\end{align*}
where $R(t)$ contains all terms dependant on $M$. Hence we get $f_{2}(t)=f(t)+S(t)$
where \[
S(t):=\begin{cases}
\frac{f(t)}{|f(t)|}\tau R(t) & t\in[0,b]\\
0 & \textrm{otherwise.}\end{cases}\]
 The next step is to ask how big does $M$ need be to ensure that
$S$ lives only in the high end of the spectrum. This is pretty straightforward,
but let us do it in detail nonetheless.

Let therefore $|n|\leq N$ and examine $\widehat{S}(n)$. Integrating
by parts we get\[
\widehat{S}(n)=\int_{0}^{1}S(t)e^{-2\pi int}=\int_{0}^{1}S(t)\, dt+2\pi in\int_{0}^{1}e^{-2\pi ins}\int_{0}^{s}S(t)\, dt\, ds.\]
To estimate this we start by writing, for $s\in[0,b]$,\begin{equation}
\left|\int_{0}^{s}S(t)\right|\leq\sum_{j=1}^{\infty}\frac{1}{j!}\sideset{}{'}\sum_{k=0}^{j}\binom{j}{k}2^{-j}\tau\left|\int_{0}^{s}e^{iM(j-2k)t}\frac{f(t)}{|f(t)|}\delta^{j}(t)\right|\label{eq:intS}\end{equation}
where $\Sigma'$ means here that for $j$ even the sum does not contain
the term $k=j/2$. Every term on the right we again estimate by integration
by parts and get \begin{equation}
\left|\int_{0}^{s}e^{iM(j-2k)t}\frac{f(t)}{|f(t)|}\delta^{j}(t)\right|\leq\frac{2}{M}\left(C+V\left(\frac{f}{|f|}\delta^{j}\right)\right).\label{eq:int<V}\end{equation}
For $j$ even we can simply estimate the variation by the maximum
of the derivative. We use Bernstein's inequality%
\footnote{In fact here it is enough to use the trivial inequality $|f'|\leq N^{2}||f||_{\infty}$.%
} to get $|f'(t)|\leq N||f||_{\infty}\leq N(1+c_{1})$ and $|(|f|)'(t)|\leq|f'(t)|\leq CN$.
Further,\[
|\left(\delta^{2j}\right)'|\leq\frac{Cj}{\tau}|f|'\leq CjN\]
and therefore \begin{equation}
V((f/|f|)\delta^{2j})\leq||(f/|f|)\delta^{2j}||_{C^{1}}\leq CjN.\label{eq:Vjeven}\end{equation}
In the case of $j$ odd we use lemma \ref{lem:C2} and get\begin{align}
V\left(\frac{f}{|f|}\delta^{2j+1}\right) & \stackrel{(\ref{eq:prod})}{\leq}CV\left(\frac{f}{|f|}\delta^{2j}\right)+CV(\delta)\leq CjN+\sqrt{\left\Vert \psi\left(1-\frac{|f|}{\tau}\right)\right\Vert _{C^{2}}}\leq\nonumber \\
 & \stackrel{(*)}{\leq}CjN+C\sqrt{\left\Vert \,|f|\,\right\Vert _{C^{2}}+\left\Vert \,|f|\,\right\Vert _{C^{1}}^{2}}\leq CjN\label{eq:Vjodd}\end{align}
where in $(*)$ we used that $\psi'\leq C$ and $\psi''\leq C$ and
in the last inequality we again used Bernstein's inequality. Inserting
this into (\ref{eq:int<V}) and that into (\ref{eq:intS}) and summing
over $k$ and $j$ gives $|\int_{0}^{s}S(t)|\leq CN/M$. Hence $|\widehat{S}(n)|\leq CnN/M$.
Summing over $n$ we get\begin{equation}
\left\Vert \sum_{n=-N}^{N}\widehat{S}(n)e^{2\pi int}\right\Vert _{\infty}\leq\frac{CN^{3}}{M}.\label{eq:smallfreq}\end{equation}

Next we need to estimate the $|\widehat{S}(n)|$ for large $n$. The
square root in the definition of $\delta$ means $S$ is not smooth,
but we shall show that $S$ is H\"older-$\frac{1}{2}$. We remind
the reader that in general such functions are have uniformly convergent
Fourier expansion. In fact, by the Dini-Lipschitz test \cite[\S 2.71]{Z68}
\begin{equation}
\Bigg\Vert\sum_{|n|>\nu}\widehat{S}(n)e^{int}\Bigg\Vert_{\infty}\leq C\nu^{-1/2}\log\nu\Vert S\Vert_{C^{1/2}}.\label{eq:VHolder}\end{equation}
Write therefore\begin{align*}
\left\Vert e^{iM(j-2k)t}\frac{f}{|f|}\delta^{j}\right\Vert _{C^{1/2}} & \stackrel{(\ref{eq:prod})}{\leq}\left\Vert \frac{f}{|f|}\delta^{j}\right\Vert _{C^{1/2}}+C\left\Vert e^{iM(j-2k)t}\right\Vert _{C^{1/2}}\leq\\
 & \stackrel{\hphantom{(\ref{eq:prod})}}{\leq}\left\Vert \frac{f}{|f|}\delta^{j}\right\Vert _{C^{1/2}}+C\sqrt{jM}\end{align*}
To estimate the term $(f/|f|)\delta^{j}$, note that for even $j$
we have \begin{equation}
\left\Vert \frac{f}{|f|}\delta^{2j}\right\Vert _{C^{1/2}}\leq\left\Vert \frac{f}{|f|}\delta^{2j}\right\Vert _{C^{1}}\stackrel{(\ref{eq:Vjeven})}{\leq}CjN.\label{eq:HolderHalf}\end{equation}
Further, since \[
\left\Vert \delta\right\Vert _{C^{1/2}}\leq\sqrt{\left\Vert \psi\left(1-\frac{|f|}{\tau}\right)\right\Vert _{1}}\leq C\]
we see that (\ref{eq:HolderHalf}) holds for $2j+1$ as well. Hence
\[
\left\Vert e^{iM(j-2k)t}\frac{f}{|f|}\delta^{j}\right\Vert _{C^{1/2}}\leq CjN+C\sqrt{jM}.\]
we now sum over $k$ and $j$ with the final result being $\left\Vert S\right\Vert _{C^{1/2}}\leq C(N+\sqrt{M})$.
Returning to (\ref{eq:VHolder}) this gives  \begin{equation}
\Bigg\Vert\sum_{|n|>\nu}\widehat{S}(n)e^{int}\Bigg\Vert_{\infty}\leq C\nu^{-1/2}\log\nu(N+\sqrt{M}).\label{eq:Holder-V}\end{equation}
Combining (\ref{eq:smallfreq}) and (\ref{eq:Holder-V}) allows us
to define\[
F=f+\sum_{N<|n|\leq M^{4}}\widehat{S}(n)e^{int}\]
and get $\left\Vert F-f_{2}\right\Vert _{\infty}\leq CN^{3}/M$. We
pick $M=\left\lceil CN^{3}/\epsilon\right\rceil $, get\begin{equation}
\left\Vert F-f_{2}\right\Vert _{\infty}\leq\epsilon\label{eq:f1-F}\end{equation}
and the lemma is proved (remember (\ref{eq:F>}) and (\ref{eq:|F|-|f|})).
\end{proof}
\begin{lem}
\label{lem:E}With the notations of lemma \ref{lem:arcsin}, if one
replaces $[a,b]$ with a simple set $E$ (satisfying that $|f|$ is
constant on $\partial E$) then the lemma holds with condition (\ref{eq:degF})
replaced by \begin{equation}
\deg F\leq C(\deg f)^{16}/\epsilon^{4}.\label{eq:specE}\end{equation}

\end{lem}
\begin{proof}
We note that $|f|^{2}$ is a trigonometric polynomial of degree $\leq2\deg f$
and hence the number of solutions of $|f|^{2}=\tau$ for any number
$\tau$ is $\leq4\deg f+1$ (we assume here that $|f|$ is not constant
--- if it is, just take $F=f$). Therefore $E$ is composed of no
more than $2\deg f$ intervals $I_{k}$. Apply lemma \ref{lem:arcsin}
for each interval $I_{k}$ with $\epsilon_{\textrm{lemma \ref{lem:arcsin}}}=\epsilon/2\deg f$.
Call the resulting functions $F(I_{k})$ and define $F=f+\sum_{k}(F(I_{k})-f)$.
All the conditions are obviously satisfied.
\end{proof}
\begin{lem}
Let $|\eta|<\tau$ be two numbers, $\eta\in\mathbb{C}$ and $\tau\in\mathbb{R}^{+}$.
Then\begin{equation}
\left|\eta\pm\sqrt{\tau^{2}-|\eta|^{2}}\frac{i\eta}{|\eta|}\right|=\tau,\label{eq:eta-tau}\end{equation}
and for any $\sigma\in[-1,1]$,\begin{equation}
|\eta|\leq\left|\eta+\sigma\sqrt{\tau^{2}-|\eta|^{2}}\frac{i\eta}{|\eta|}\right|\leq\tau.\label{eq:eta<tau}\end{equation}

\end{lem}
This follows from Pythagoras' theorem since we are adding orthogonal
vectors.

\begin{lem}
\label{lem:smooth}Let $f$, $\epsilon$ and $[a,b]$ be as in lemma
\ref{lem:arcsin}, and let $[a',b']$ be another interval, $b'-a'\geq b-a$.
Let $g$ be a polynomial satisfying $|\widehat{f}(n)|=|\widehat{g}(n)|$
for all $n$. Then one can extend $f$ and $g$ as $F$ and $G$ still
satisfying $|\widehat{F}(n)|=|\widehat{G}(n)|$ and such that (\ref{eq:lem|F|-|f|}),
(\ref{eq:F-f}) hold as well as\begin{equation}
\begin{aligned}|g(t)-G(t)| & <C\sqrt{(b-a)/(b'-a')} & \quad & \forall t\in[a',b']\\
|g(t)-G(t)| & <\epsilon &  & \forall t\not\in[a',b'].\end{aligned}
\label{eq:g-G}\end{equation}

\end{lem}
\begin{proof}
Assume w.l.o.g.~that $a=a'=0$. Let $l$ satisfy that $4^{-l}\geq b/b'>4^{-l-1}$,
and assume $l>0$ (otherwise one can take $F$ from lemma \ref{lem:arcsin}
and $G=g+F-f$). Let $\delta\in\left]0,\epsilon\right[$ and $M$
be some parameters to be fixed later --- $\delta$ will be taken sufficiently
small and $M$ sufficiently large, depending on $\delta$. Let $\psi$
be an $M$-approximation of the first Radamacher function $r_{1}:=\mathbf{1}_{[0,1/2]}-\mathbf{1}_{[1/2,1]}$
namely\[
\psi=\sum_{|n|\leq M}\widehat{r_{1}}(n)\frac{M-|n|}{M}e^{int}.\]
As is well known, $\psi$ is real, $|\psi|\leq1$ and \begin{equation}
|r_{1}(t)-\psi(t)|\leq CM^{-1/2}\quad\forall t,\,\left\langle 2t\right\rangle >M^{-1/2}\label{eq:r1-psi}\end{equation}
 where $\left\langle x\right\rangle $ is defined (somewhat non-standardly)
as the distance of $x$ from the integers. For any integer $m$ define
$\psi_{[m]}$ using\[
\psi_{[m]}(t)=\psi(mt)\]
where we understand here $\psi$ as a $1$-periodic function. We define
$3l$ functions which mimic the behavior of $3l$ different Radamacher
functions:\[
s_{i}:=\psi_{[(3M)^{i}]}\quad i=1,\dotsc,3l.\]
The reason we are approximating the Radamacher functions is the following
innocuous equality\[
|r_{1}+r_{2}+r_{3}-r_{1}r_{2}r_{3}|=2\]
which holds at every point except the jump points. This of course
holds for other triplets i.e. $|r_{3i-2}+r_{3i-1}+r_{3i}-r_{3i-2}r_{3i-1}r_{3i}|=2$.
With this in mind let us construct $4^{l}$ functions which mimic
products of the four Walsh functions $r_{1},r_{2},r_{3},-r_{1}r_{2}r_{3}$.
Formally, for every sequence $\{\epsilon_{i}\}_{i=1}^{l}$with $\epsilon_{i}\in\{1,2,3,4\}$
we define\[
\sigma_{\{\epsilon_{i}\}}=2^{-l}\prod_{i=1}^{l}\begin{cases}
s_{3(i-1)+\epsilon_{i}} & \epsilon_{i}=1,2,3\\
-s_{3i-2}s_{3i-1}s_{3i} & \epsilon_{i}=4.\end{cases}\]
For convenience, if $j=\sum_{i=1}^{\alpha}(\epsilon_{i}-1)4^{i-1}$
is some number between $0$ and $4^{l}-1$ we will denote $\sigma_{j}:=\sigma_{\{\epsilon_{i}\}}$.
The $\sigma_{j}$-s satisfy the following properties:
\begin{enumerate}
\item \label{enu:sigma-2^-k}For all $j$ and $t$,\begin{equation}
|\sigma_{j}(t)|\leq2^{-l}\label{eq:sigma-2^-l}\end{equation}

\item Let $\mathcal{B}$ be the set of {}``bad'' $t$-s satisfying $\left\langle 2(3M)^{i}t\right\rangle \leq M^{-1/2}$
for some $i=1,\dotsc,3l$. Then \begin{equation}
1-ClM^{-1/2}<\bigg|\sum_{j}\sigma_{j}(t)\bigg|\leq1\quad\forall t\not\in\mathcal{B}.\label{eq:sumj}\end{equation}
This follows because\[
\sum_{j}\sigma_{j}(t)=2^{-l}\prod_{i=1}^{l}(s_{3i-2}+s_{3i-1}+s_{3i}-s_{3i-2}s_{3i-1}s_{3i}).\]
As explained above, $|r_{1}+r_{2}+r_{3}-r_{1}r_{2}r_{3}|=2$ and hence
so does each term in the product (with an error of $CM^{-1/2})$ and
the product has absolute value $2^{l}$. Further, an easy calculation
shows that $x+y+z-xyz\leq2$ for all $x,y,z\in[-1,1]$ so \begin{equation}
\bigg|\sum_{j}\sigma_{j}(t)\bigg|\leq1\quad\forall t\label{eq:sumsig1}\end{equation}
 $t$ both outside and inside $\mathcal{B}$.
\item \label{enu:disjoint}The spectra of $\sigma_{j}$ are disjoint. If
$M>3\deg f$ then they are also disjoint from $\spec f$. This is
also easy --- in fact the spectra of any product of $s_{i}$-s are
disjoint, and disjoint from $\spec f$.
\end{enumerate}
Let $\tau:=|f(0)|$, let $\varphi$ be defined by \[
\varphi(t)=\begin{cases}
\sqrt{\tau^{2}-|f(t)|^{2}}\frac{if(t)}{|f(t)|} & t\in[0,b]\\
0 & \textrm{otherwise,}\end{cases}\]
and let $P$ be a trigonometric polynomial approximating $\varphi$,
$||P-\varphi||_{\infty}<\delta$. Assume $M>3\deg P$ so that the
spectra of $P\sigma_{j}$ are all disjoint and disjoint from $\spec f$. 

We are now in a position to define our first approximation step, \begin{align*}
f_{2}(t) & :=f(t)+\sum_{j=0}^{4^{l}-1}P\sigma_{j}(t)\\
g_{2}(t) & :=g(t)+\sum_{j=0}^{4^{l}-1}(P\sigma_{j})(t-jb).\end{align*}
It is clear that $|\widehat{f_{2}}(n)|=|\widehat{g_{2}}(n)|$ since
the only difference between them is a rotation of each $P\sigma_{j}$.
Since they have disjoint spectra, this preserves the absolute value
of the Fourier transform.

Examine one $t\in[0,b]\setminus\mathcal{B}$. We use (\ref{eq:sumj})
to sum over the $j$-s and get\[
f_{2}=f\pm\varphi(t)+R(t)\quad|R|\leq ClM^{-1/2}+\delta\quad\forall t\in[0,b]\setminus\mathcal{B}.\]
Notice that by Pythagoras' theorem (\ref{eq:eta-tau}) for every $t\in[0,b]$
$|f\pm\varphi|=\tau$, so \begin{equation}
|\,|f_{2}(t)|-\tau|\leq ClM^{-1/2}+\delta\quad\forall t\in[0,b]\setminus\mathcal{B}.\label{eq:f2 f pm phi}\end{equation}
On $[0,b]\cap\mathcal{B}$ we use (\ref{eq:eta<tau}) and (\ref{eq:sumsig1})
to get \begin{equation}
|f(t)|-\delta<|f_{2}(t)|<\tau+\delta.\label{eq:f2f on B}\end{equation}
Finally, outside $[0,b]$ we have $|f_{2}(t)-f(t)|\leq\delta$ regardless
of whether $t\in\mathcal{B}$ or not.

As for $g_{2}$, because the various translates $\varphi(t-jb)$ have
disjoint support we get (remember (\ref{eq:sigma-2^-l}))\begin{equation}
\begin{aligned}|g_{2}(t)-g(t)| & \leq2^{-l}\max\varphi+2^{l}\delta\leq C\sqrt{b/b'}+2^{l}\delta & \quad & \forall t\in[0,b']\\
|g_{2}(t)-g(t)| & \leq2^{l}\delta &  & \forall t\not\in[0,b'].\end{aligned}
\label{eq:g2-g}\end{equation}
These are the properties we need for $f_{2}$ and $g_{2}$.

\smallskip{}

\noindent \textbf{Step 2}: $f_{2}$ and $g_{2}$ satisfy the conditions
of the lemma except on the small set $\mathcal{B}$. On it we correct
using lemma \ref{lem:E}. Assume therefore that $\delta$ is sufficiently
small and $M$ sufficiently large such that \begin{equation}
1-c_{1}<|f_{2}|<1+c_{1}\label{eq:f2<1+c1}\end{equation}
 and we can use lemma \ref{lem:E} for $f_{2}$. Apply it with the
parameter $\epsilon_{\textrm{lemma \ref{lem:E}}}=\frac{1}{2}\epsilon$
and the set \[
E=\{ t\in[0,b]\cap\mathcal{B}:|f_{2}(t)|\leq\tau-{\textstyle \frac{1}{2}}\epsilon\}.\]
If $\delta$ is sufficiently small and $M$ sufficiently large then
by (\ref{eq:f2 f pm phi}) we have that $E$ is contained in the interior
of $\mathcal{B}$ and therefore the condition that $|f_{2}|$ is constant
on $\partial E$ is satisfied. Call the resulting function $F$. We
get\begin{alignat}{2}
|\,|F(t)|-\tau| & \leq\epsilon & \quad & \forall t\in E\label{eq:|F|tau E}\\
|F(t)-f_{2}(t)| & \leq{\textstyle \frac{1}{2}}\epsilon &  & \forall t\not\in E.\label{eq:F-f off E}\end{alignat}
As for $\deg F$, since $\deg f_{2}\leq(3M)^{3l+1}$ we get \[
\deg F\leq C(3M)^{48l+16}/\epsilon^{4}.\]
This gives all required properties from $F$. Hence we need to define
$G$. For this purpose examine the random function\[
h:=\sum_{n}\pm\widehat{F-f_{2}}(n)e^{int}\]
where the $\pm$ stands for independent Bernoulli variables. It is
well known (see \cite[chapter 6, theorem 2]{K85}) that with high
probability \[
\left\Vert h\right\Vert _{\infty}\leq C\left\Vert F-f_{2}\right\Vert _{2}\sqrt{\log\deg F}\]
and in particular there exists a choice of signs $\xi_{n}$ satisfying
this inequality. Define \[
G:=g_{2}+\sum_{n}\xi_{n}\widehat{F-f_{2}}(n)e^{2\pi int}.\]
Clearly we have $|\widehat{F}(n)|=|\widehat{G}(n)|$. Further, $\left\Vert F-f_{2}\right\Vert _{2}\leq C\sqrt{|\mathcal{B}|}\leq Cl^{1/2}M^{-1/4}$
so\begin{equation}
\left\Vert G-g_{2}\right\Vert _{\infty}\leq Cl^{1/2}M^{-1/4}\sqrt{l\log M+\log1/\epsilon}.\label{eq:ClM}\end{equation}
All that is required is to pick $\delta$ and $M$ correctly. Requirements
(\ref{eq:lem|F|-|f|}) and (\ref{eq:F-f}) will follow if only $ClM^{-1/2}+\delta+\frac{1}{2}\epsilon<\epsilon$.
To see (\ref{eq:lem|F|-|f|}) note that on $[0,b]\setminus\mathcal{B}$
it follows from (\ref{eq:f2 f pm phi}) and (\ref{eq:F-f off E}).
On $[0,b]\cap E$ it follows from (\ref{eq:|F|tau E}). And on $[0,b]\cap(\mathcal{B}\setminus E)$
from (\ref{eq:f2f on B}), (\ref{eq:F-f off E}) and the definition
of $E$. Seeing requirement (\ref{eq:F-f}) is similar. Next, (\ref{eq:g-G})
will follow if the right hand side of (\ref{eq:ClM}) is $<\min\{\frac{1}{2}\epsilon,\sqrt{b/b'}\}$
and $2^{l}\delta<\min\{\frac{1}{2}\epsilon,\sqrt{b/b'}\}$ (remember
(\ref{eq:g2-g})). We remind the reader that in addition we assumed
$M>3\deg f$, $M>3\deg P$, that \[
\delta+ClM^{-1/2}<c_{1}-\max_{t}|\,|f(t)|-1|\]
 which ensures (\ref{eq:f2<1+c1}), and that $\delta+ClM^{-1/2}<\frac{1}{2}\epsilon$
which ensures that $|f_{2}|$ is constant on $\partial E$. Clearly
choosing $\delta$ sufficiently small and then $M$ sufficiently large
depending on $\delta$ (the most important dependency is via $M>3\deg P$)
will satisfy all these requirements and prove the lemma.
\end{proof}

\subsection{Intermediate remarks }

Since you reached so far down in the proof itself we believe you might
be interested in some remarks on the structure of the proof more substantial
than the ones given in the introduction. We start with a remark on
lemma \ref{lem:smooth}. A tempting simplification is as follows:
use lemma \ref{lem:arcsin} to correct $f$ to $F$ and then the probabilistic
argument above to construct $G=g+\sum\pm\widehat{F-f}(n)e^{int}$.
This would give the lemma with an additional benign-looking factor
of $\sqrt{\log\deg F}$. However, this $\sqrt{\log}$ factor is not
so easy to get rid of! To appreciate how serious a burden was removed,
try to estimate the relation between $n=\deg f$ and $N=\deg F$ in
lemma \ref{lem:killf} below. We got the tetration\[
N\approx\underbrace{n^{n^{n^{.^{.^{.^{n}}}}}}}_{C\log^{2}1/\epsilon\textrm{ times}}\]
(this is after some optimizations, directly following the proof would
give much more). This, by the way, is also the best we can say about
the smoothness of the final $F$ and $G$ i.e.~$\omega(\delta;F)$
and $\omega(\delta;G)$ decrease like an inverse tetration.

This is why we chose the current approach, and starting from lemma
\ref{lem:nologo} we no longer need to control the spectrum of $F$.
Put differently, $\deg F$ is the only parameter which gets worse
when one increases the various parameters of our construction (e.g.
the $l$ and $M$ of lemma \ref{lem:smooth}, the $l$, $N$ and $M$
of lemma \ref{lem:walshIJ} below etc). Removing the requirement to
control $\deg F$ gives us the flexibility to increase these parameters
with no punishment.

\subsection{The global correction scheme}

\begin{lem}
\label{lem:nologo}Lemma \ref{lem:smooth} holds with $[a,b]$ replaced
by a general set $E$ if (\ref{eq:g-G}) is replaced by\[
||g-G||_{\infty}\leq C\sqrt{|E|}.\]

\end{lem}
Note that there is no $[a',b']$ in this formulation (or rather it
is $[0,1]$).

\begin{proof}
Write $E$ as a disjoint union $E=I_{1}\cup\dotsb\cup I_{N}$ and
let $J_{1},\dotsc,J_{N}$ be disjoint intervals with $|J_{i}|=|I_{i}|/|E|$
(so they cover $[0,1]$). Denote by $\tau$ the common value of $|f(t)|$
for all $t\in\partial E$. Let $\epsilon_{2}$ be sufficiently small
such that for any $\delta\leq\epsilon_{2}$ and any interval $I\subset E$,
$\{ t\in I:|f(t)|<\tau-\delta\}$ is an interval. Let $\epsilon_{3}=\min\{\epsilon,\epsilon_{2},\sqrt{|E|}\}$.
Now use lemma \ref{lem:smooth} inductively $N$ times to get functions
$f_{i}$, $g_{i}$ satisfying
\begin{enumerate}
\item \label{enu:inductstart}$|\,|f_{i}(t)|-\tau\,|<\epsilon_{3}(1-2^{-i-2})$
for all $t\in\bigcup_{j=1}^{i}I_{j}$,
\item $|f_{i}(t)-f(t)|\leq\epsilon_{3}(\frac{1}{4}-2^{-i-2})$ for all other
$t$-s,
\item $|g_{i}(t)-g(t)|\leq C\sqrt{|E|}+\epsilon_{3}(\frac{1}{4}-2^{-i-2})$
for all $t\in\bigcup_{j=1}^{i}J_{j}$, and
\item \label{enu:inductend}$|g_{i}(t)-g(t)|\leq\epsilon_{3}(\frac{1}{4}-2^{-i-2})$
for all other $t$-s.
\end{enumerate}
\noindent Initiate the induction with $f_{0}=f$, $g_{0}=g$. We need
to define the parameters for lemma \ref{lem:smooth}, most importantly
the interval. Examine therefore $I_{i}$. Use the definition of $\epsilon_{2}$
with $\delta=\frac{3}{4}\epsilon_{3}$ to get that \[
I_{i}^{*}:=\{ t\in I_{i}:|f(t)|<\tau-{\textstyle \frac{3}{4}}\epsilon_{3}\}\]
 is an interval. By the second induction assumption we know that $|f_{i-1}(t)|<\tau-\frac{1}{2}\epsilon_{3}$
for all $t\in I_{i}^{*}$ and $|f_{i-1}(t)|>\tau-\epsilon_{3}(1-2^{-i-1})$
for all $t\in I_{i}\setminus I_{i}^{*}$. Let $I_{i}^{**}$ be the
component of $\{ t\in I_{i}:|f_{i-1}(t)|\leq\tau-\frac{1}{2}\epsilon_{3}\}$
containing $I_{i}^{*}$. Now use lemma \ref{lem:smooth} with the
parameters as in the following table:%
\begin{table}[H]
\begin{tabular}{|c||c|c|c|c|c||c|c|}
\hline 
lemma \ref{lem:smooth}&
$f$&
$g$&
$\epsilon$&
$[a,b]$&
$[a',b']$&
$F$&
$G$\tabularnewline
\hline 
here&
$f_{i-1}$&
$g_{i-1}$&
$\epsilon_{3}2^{-i-2}$&
$I_{i}^{**}$&
$J_{i}$&
$f_{i}$&
$g_{i}$\tabularnewline
\hline
\end{tabular}
\end{table}

\noindent i.e.~the lemma's output will be used to define $f_{i}$
and $g_{i}$. It is easy to verify \ref{enu:inductstart}-\ref{enu:inductend}
and the induction is complete. $F=f_{N}$ and $G=g_{N}$ are the desired
functions.
\end{proof}
\begin{defn*}
Let $f$ be a trigonometric polynomial, let $E$ be a set and let
$\max_{E}|f|\leq\alpha<\beta$. Let $\epsilon>0$. Then an \textbf{$\epsilon$-lifting
of $f$ on $E$ from $\alpha$ to $\beta$} is a trigonometric polynomial
$F$ extending $f$ such that
\begin{enumerate}
\item \label{enu:|f-F| not E}$|f(t)-F(t)|<\epsilon$ for all $t\not\in E$,
\item \label{enu:E minus eps}For all $t\in E$ except a set of measure
$<\epsilon$,\begin{equation}
(\beta-\alpha)+|f(t)|-\epsilon<|F(t)|<\beta+\epsilon,\label{eq:lifting}\end{equation}

\item \label{enu:IJ F-f all}For all $t\in E$, \[
\hphantom{(\beta-\alpha)\:+\:}|f(t)|-\epsilon<|F(t)|<\beta+\epsilon,\]

\item \label{enu:||F-f||}$\left\Vert F-f\right\Vert _{\infty}<2\sqrt{\beta^{2}-\alpha^{2}}$.
\end{enumerate}
See figure 1. %
\begin{figure}
\input{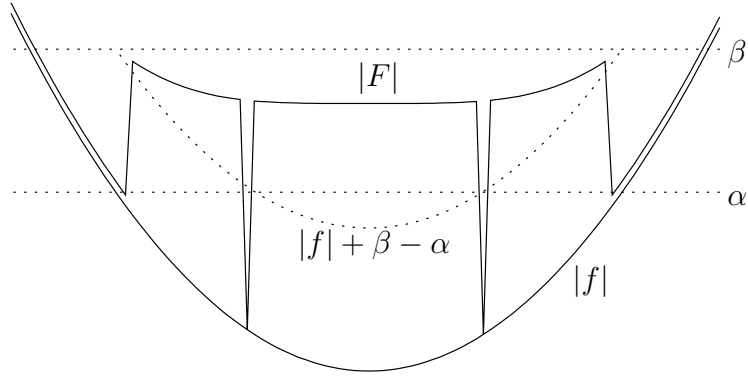}

\caption{$F$ is a lifting of $f$ from $\alpha$ to $\beta$ on the set $E$
where $|f|\leq\alpha$. The graph depicts only the absolute values
of $f$ and $F$. $|F|$ is drawn here taking the value $\sqrt{|f|^{2}+\beta^{2}-\alpha^{2}}$
which is what will typically happen in our construction.}
\end{figure}
If $F$ and $G$ are two liftings of $f$ and $g$ which satisfy $|\widehat{F}(n)|=|\widehat{G}(n)|$
for all $n$ then we call them {}``compatible liftings''.
\end{defn*}
\begin{lem}
\label{lem:walshIJ}Let $f$ and $g$ be as in lemma \ref{lem:smooth}.
Let $\epsilon>0$ be some parameter. Let $l\geq0$ be some integer,
let $I$ be an interval with length $\leq4^{-l}$ and let $J$ be
an interval of length $4^{l}|I|$. Assume that $|f(t)|<\tau$ for
all $t\in I$ where $\tau<1$. Then there exist compatible $\epsilon$-liftings
$F$ of $f$ on $I$ from $\tau$ to $1$ and $G$ of $g$ on $J$
from any $\nu\geq\max_{J}|g|$ to $\mu=\sqrt{\nu^{2}+4^{-l}(1-\tau^{2})}$.
\end{lem}
\begin{proof}
We may assume $I=[0,b]$ and $J=[0,4^{l}b]$. Define two parameters
$N$ and $M$ that will accompany us throughout the proof and will
be fixed at the end (think of both as large but of $M$ as being much
larger than $N$). The proof is very similar to the proof of lemma
\ref{lem:smooth}, and in particular we retain the notations $\psi$
and $\psi_{[m]}$. We shall repeat the construction of the $\sigma_{j}$-s
of lemma \ref{lem:smooth}, $N$ times in disjoint spectra. Namely,
we define functions $\sigma_{j}^{q}$, $j=0,\dotsc,4^{l}-1$ and $q=0,\dotsc,N-1$
satisfying
\begin{enumerate}
\item Let $\mathcal{B}$ be the set of {}``bad'' $t$-s satisfying $\left\langle 2(3M)^{i}t\right\rangle \leq M^{-1/2}$
for some $i=1,\dotsc,3lN$. Then \begin{equation}
2^{-l}-ClM^{-1/2}<|\sigma_{j}^{q}(t)|\leq2^{-l}\quad\forall j,\, q,\, t\not\in\mathcal{B}.\label{eq:sigma jq-2^-l}\end{equation}
 Further, $|\sigma_{j}^{q}(t)|\leq2^{-l}$ for all $t$.
\item For any $t$ and $q$, $|\sum_{j}\sigma_{j}^{q}(t)|\leq1$ and further\begin{equation}
1-ClM^{-1/2}<\bigg|\sum_{j}\sigma_{j}^{q}(t)\bigg|\leq1\quad\forall q,\, t\not\in\mathcal{B}.\label{eq:sumjq}\end{equation}

\item The spectra of $P\sigma_{j}^{q}$ are disjoint. Here $P$ is any polynomial
with $\deg P\leq M$.
\end{enumerate}
Let $P$ be an $M$-approximation of the indicator function $\mathbf{1}_{[0,b/N]}$
(in the same sense of $\psi$ in lemma \ref{lem:smooth}) so that
$P$ is real, $|P|\leq1$ and \begin{equation}
|P(t)-\mathbf{1}_{[0,b/N]}(t)|\leq CM^{-1/2}+C\frac{N/b}{M}\quad\forall t,d(t,\{0,b/N\})>M^{-1/2}\label{eq:phi-1}\end{equation}
where $d$ is the usual distance of a point from a set considered
periodically. We will assume henceforth that $M>N^{2}b^{-2}$ and
avoid carrying the $N/bM$ term. This allows us to define our functions\begin{align}
F & =f+\sqrt{1-\tau^{2}}\sum_{q=0}^{N-1}\sum_{j=0}^{4^{l}-1}\frac{if(qb/N)}{|f(qb/N)|}(P\sigma_{j}^{q})(t-qb/N)\label{eq:def F}\\
G & =g+\sqrt{1-\tau^{2}}\sum_{q=0}^{N-1}\sum_{j=0}^{4^{l}-1}\frac{ig(jb+qb/N)}{|g(jb+qb/N)|}(P\sigma_{j}^{q})(t-(jb+qb/N)).\label{eq:def G}\end{align}
As in lemma \ref{lem:smooth} it is clear that $|\widehat{F}(n)|=|\widehat{G}(n)|$
since each term is rotated and multiplied by a unimodular number.
Define \[
\mathcal{B}':=\bigg(\bigcup_{q=0}^{N-1}\bigcup_{j=0}^{4^{l}-1}(\mathcal{B}\cup[-M^{-1/2},M^{-1/2}])+jb+qb/N\bigg)\]
where the $+$ here stands for usual set addition i.e.~$B+x:=\{ b+x:b\in B\}$.
Examine some $t\not\in\mathcal{B}'$. We use (\ref{eq:sumjq}) to
sum over the $j$-s in (\ref{eq:def F}) and get\begin{equation}
F=f\pm\sqrt{1-\tau^{2}}\sum_{q=0}^{N-1}\frac{if(qb/N)}{|f(qb/N)|}P(t-qb/N)+O(lNM^{-1/2}).\label{eq:linerror}\end{equation}
By (\ref{eq:phi-1}) all the terms in the sum over $q$ give $P=O(M^{-1/2})$
unless $t\in[0,b]$ in which case the term $q=\left\lfloor Nt/b\right\rfloor $
gives $P=1+O(M^{-1/2})$. We get \begin{equation}
F(t)=f(t)\pm\sqrt{1-\tau^{2}}\frac{if(qb/N)}{|f(qb/N)|}+O(lNM^{-1/2})\quad\forall t\in[0,b]\setminus\mathcal{B}'.\label{eq:Ffpmsqrt}\end{equation}
Now we use Pythagoras' theorem (\ref{eq:eta-tau}) for $\eta=\tau f(qb/N)/|f(qb/N)|$
and get \begin{multline}
\left|1-\left|f(qb/N)\pm\sqrt{1-|\eta|^{2}}\frac{i\eta}{|\eta|}\right|\,\right|\leq\\
\left|\left(\eta\pm\sqrt{1-|\eta|^{2}}\frac{i\eta}{|\eta|}\right)-\left(f(qb/N)\pm\sqrt{1-|\eta|^{2}}\frac{i\eta}{|\eta|}\right)\right|=\tau-|f(qb/N)|\label{eq:fqbN}\end{multline}
which allows to estimate, for any $t\in[0,b]\setminus\mathcal{B}'$,
\begin{align*}
|F(t)| & \stackrel{(\ref{eq:Ffpmsqrt})}{\geq}\left|f\left(\frac{qb}{N}\right)\pm\sqrt{1-|\eta|^{2}}\frac{i\eta}{|\eta|}\right|-\left|f(t)-f\left(\frac{qb}{N}\right)\right|-ClNM^{-1/2}\geq\\
 & \stackrel{(\ref{eq:fqbN})}{\geq}1-\tau+\left|f\left(\frac{qb}{N}\right)\right|-\left|f(t)-f\left(\frac{qb}{N}\right)\right|-ClNM^{-1/2}\geq\\
 & \stackrel{\hphantom{(\ref{eq:fqbN})}}{\geq}1-\tau+|f(t)|-2\omega(b/N;f)-ClNM^{-1/2}.\end{align*}
Denote $\omega:=2\omega(b/N;f)+ClNM^{-1/2}$. In the other direction
(\ref{eq:eta<tau}) gives $|f(qb/N)|\leq\left|f(qb/N)+\sigma\sqrt{1-\tau^{2}}\frac{i\eta}{|\eta|}\right|\leq1$
for any $\sigma\in[-1,1]$. We take the variation in $f$ into consideration
as above and get\begin{alignat*}{2}
(1-\tau)+|f(t)|-\omega & \leq|F(t)|\leq1+\omega & \quad & \forall t\in[0,b]\setminus\mathcal{B}'\\
|f(t)|-\omega & \leq|F(t)|\leq1+\omega & \quad & \forall t\in[0,b]\cap\mathcal{B}'.\end{alignat*}
Finally we also have $\left\Vert F-f\right\Vert _{\infty}\leq\sqrt{1-\tau^{2}}+\omega$.

Next we move to examine $G$. This time we first notice that the only
meaningful term in (\ref{eq:def G}) is the one for which $j=\left\lfloor t/b\right\rfloor $
and $q=\left\lfloor N(t/b-j)\right\rfloor $. So we have (again for
$t\not\in\mathcal{B}'$)\[
G(t)=g(t)+\sqrt{1-\tau^{2}}\frac{ig(t+R(t))}{|g(t+R(t))|}\sigma_{\left\lfloor t/b\right\rfloor }(R(t))+O(2^{l}NM^{-1/2})\quad|R(t)|\leq b/N\]
(we get here a $2^{l}$ factor in the error, compared to the $l$
in (\ref{eq:linerror}) because the $\sigma_{j}$-s are no longer
{}``synchronized'' so we cannot use (\ref{eq:sumjq}) and have to
use (\ref{eq:sigma jq-2^-l}) and sum over all the terms). We now
use (\ref{eq:sigma jq-2^-l}) and get\[
G(t)=g(t)\pm2^{-l}\sqrt{1-\tau^{2}}\frac{ig(t+R(t))}{|g(t+R(t))|}+O(2^{l}NM^{-1/2})\quad\forall t\in[0,b']\setminus\mathcal{B}'.\]
A calculation similar to the one done with $f$ shows that\begin{gather*}
\begin{aligned}\mu-\nu+|g(t)|-\omega' & \leq|G(t)|\leq\mu+\omega' &  & \forall t\in[0,b']\setminus\mathcal{B}'\\
|g(t)|-\omega' & \leq|G(t)|\leq\mu+\omega' &  & \forall t\in[0,b']\cap\mathcal{B}'\end{aligned}
\\
\omega':=2\omega(b/N;g)+C2^{l}NM^{-1/2}\end{gather*}
and $\left\Vert G-g\right\Vert _{\infty}\leq\sqrt{1-\tau^{2}}+\omega'$. 

With these estimates the lemma will be finished once we pick $N$
and $M$. First pick $N$ such that $\omega(b/N;f)<\min\{\frac{1}{4}\epsilon,\frac{1}{4}\sqrt{1-\tau^{2}}\}$
and similarly for $g$. Next pick $M$ to satisfy all past requirements.
They are all of the type {}``$M$ is sufficiently large (possibly
depending on $N$, $l$ and $\epsilon$)''. Here is the full list
(in chronological order): $M>N^{2}b^{-2}$, $M>3\deg f$, $ClNM^{-1/2}<\min\{\frac{1}{2}\epsilon,\frac{1}{2}\sqrt{1-\tau^{2}}\}$
(which bounds $\omega$ and ensures $|F(t)-f(t)|<\epsilon$ outside
$[0,b]$) and $C4^{l}NM^{-1/2}<\epsilon$ which ensures $|\mathcal{B}'|<\epsilon$,
$\omega'<\epsilon$ and $|G(t)-g(t)|<\epsilon$ outside $[0,b']$.
With all these satisfied we get everything we want for $F$ and $G$.
\end{proof}
Again we need to generate a set version from the interval version.
We trust that by now the reader will have no problem to prove:

\begin{lem}
\label{lem:walshE}Let $f$, $g$, $\epsilon$, $\tau$ and $l$ be
as in lemma \ref{lem:walshIJ}. Let $E$ be a simple set with $|E|=4^{-l}$
such that $|f(t)|<\tau$ for all $t\in E$. Then one can find compatible
$\epsilon$-liftings $F$ of $f$ on $E$ from $\tau$ to $1$ and
$G$ of $g$ on $[0,1]$ from $\nu\geq||g||_{\infty}$ to $\mu=\sqrt{\nu^{2}+4^{-l}(1-\tau^{2})}$.
\end{lem}
In the next lemma we get rid of the errors in the exceptional sets
(compare clauses \ref{enu:half f} and \ref{enu:G mu g nu} below
to requirement (\ref{eq:lifting}) from the definition of lifting).
Hence it will be convenient to use the following definition: The oscillation
of the absolute value of a function $g$ is\[
\Osc(g):=\max_{t}|g(t)|-\min_{t}|g(t)|\]

\begin{lem}
\label{lem:halff}Let $f$ and $g$ be as in lemma \ref{lem:smooth}
but with the additional requirements\begin{equation}
1-c_{2}<|f|<1\quad1-c_{2}<|g|<1+c_{2}\label{eq:reqc2}\end{equation}
for some absolute constant $c_{2}>0$. Let $\epsilon>0$ be some parameter.
Then one can extend $f$ and $g$ to $F$ and $G$ such that
\begin{enumerate}
\item \label{enu:half f}$1-\frac{1}{2}\left\Vert 1-|f|\,\right\Vert _{\infty}-\epsilon<|F(t)|<1+\epsilon$
for all $t$,
\item $\left\Vert F-f\right\Vert _{\infty}<C\sqrt{\left\Vert 1-|f|\,\right\Vert _{\infty}}$
\item \label{enu:G mu g nu}$\Osc(G)<\Osc(g)+\epsilon$
\item $\left\Vert G-g\right\Vert _{\infty}<C\sqrt{\left\Vert 1-|f|\,\right\Vert _{\infty}}$.
\end{enumerate}
\end{lem}

\begin{proof}
Denote $\tau=1-\frac{1}{2}\left\Vert 1-|f|\,\right\Vert _{\infty}$.
Let \[
E=\Big\{ t:|f(t)|<\tau-{\textstyle \frac{1}{4}}\epsilon\Big\}\]
and write $|E|$ in base $4$ namely\[
|E|=\sum_{l=1}^{\infty}\alpha_{l}4^{-l}\quad\alpha_{l}\in\{0,1,2,3\}.\]
Divide $E$ into simple sets $E_{l,i}$, $l=1,2,\dotsc$ and $i=0,\dotsc,\alpha_{l}-1$
with $|E_{l,i}|=4^{-l}$. Order these sets by increasing $l$ and
call the resulting sequence $\{ E_{j}\}_{j=1}^{\infty}$. We now create
two sequences of polynomials $f_{i}$ and $g_{i}$ with $|\widehat{f_{i}}(n)|=|\widehat{g_{i}}(n)|$
for all $n$ and all $i$ inductively by using lemma \ref{lem:walshE}
(at the $i$'th step) with $f_{i-1}$, $g_{i-1}$ the set $E_{i}$,
the parameters $\epsilon_{\textrm{\mbox{lemma \ref{lem:walshE}}}}=\delta2^{-i}$
($\delta$ is some parameter $<\frac{1}{7}\epsilon$ to be fixed later)
and $\tau$ and with $\nu_{\textrm{lemma \ref{lem:walshE}}}=\mu_{i-1}+\delta2^{-i}$
(the induction is initialized with $f_{0}=f$, $g_{0}=g$ and $\mu_{0}=\Vert g\Vert_{\infty}$).
Call the output of the lemma $f_{i}$, $g_{i}$ and $\mu_{i}$. It
is easy to verify that 
\begin{enumerate}
\item $\tau-\frac{1}{7}\epsilon<|f_{i}(t)|<1+\frac{1}{7}\epsilon$ except
on a set $\mathcal{B}$ of measure $\delta+4^{1-i/4}$ (the $\delta$
error is the combined error from the previous stages while $4^{1-i/4}$
is the set of $E_{j}$-s not yet handled).
\item Uniformly we have\begin{equation}
\left\Vert f_{i}-f\right\Vert _{\infty}<2\sqrt{1-\tau^{2}}+\delta.\label{eq:fi-f}\end{equation}
This requires both \ref{enu:|f-F| not E} and \ref{enu:||F-f||} from
the definition of lifting.
\item \label{enu:gimui}Outside a set of measure $<\delta$ we have\begin{align*}
|g_{i}(t)| & >|g_{i-1}(t)|+\mu_{i}-(\mu_{i-1}+\delta2^{-i})-\delta2^{-i}>\\
 & >|g_{i-2}(t)|+\mu_{i}-\mu_{i-2}-2\delta(2^{1-i}+2^{-i})>\dotsb\\
 & >|g(t)|+\mu_{i}-\Vert g\Vert_{\infty}-2\delta.\end{align*}
We prefer to write this as\begin{equation}
\mu_{i}-\Osc(g)-\frac{2}{7}\epsilon<|g_{i}(t)|<\mu_{i}+\frac{1}{7}\epsilon\label{eq:giosc}\end{equation}

\item Uniformly\begin{equation}
\left\Vert g_{i}-g\right\Vert _{\infty}<2\sqrt{1-\tau^{2}}+\delta.\label{eq:gi-g}\end{equation}

\end{enumerate}
Take some $i$ sufficiently large (to be fixed later) and examine
$f_{i}$ and $g_{i}$.

\textbf{Step 2}: We now correct over the exceptional sets using lemma
\ref{lem:nologo} twice. First use the lemma with $f_{i}$, $g_{i}$,
the set $\mathcal{B}:=\{ t:|f_{i}(t)|\leq\tau-\frac{1}{7}\epsilon\}$
and $\epsilon_{\textrm{lemma }\ref{lem:nologo}}=\frac{1}{7}\epsilon$.
To enable this, fix \begin{equation}
c_{2}:=\min\left\{ \frac{1}{4}c_{1},1-\sqrt{1-\frac{c_{1}^{2}}{64}}\right\} .\label{eq:defc2}\end{equation}
With this value, the fact that $\tau>1-c_{2}$ implies that $2\sqrt{1-\tau^{2}}<\frac{1}{4}c_{1}$
and hence for $\delta<\frac{1}{4}c_{1}$ (\ref{eq:fi-f}) and (\ref{eq:reqc2})
give $1-c_{1}<|f_{i}|<1+c_{1}$ and lemma \ref{lem:nologo} may indeed
be applied. Call the resulting functions $F^{*}$ and $G^{*}$. Since
$|\mathcal{B}|<\delta+4^{1-i/4}$ it is clear that for $\delta$ sufficiently
small and $i$ sufficiently large we get \begin{equation}
\left\Vert g_{i}-G^{*}\right\Vert _{\infty}<\min\{{\textstyle \frac{1}{7}}\epsilon,{\textstyle \frac{1}{4}}c_{1}\}.\label{eq:gG*}\end{equation}
This corrects $F^{*}$ in the sense that now \begin{equation}
\tau-{\textstyle \frac{2}{7}}\epsilon<|F^{*}(t)|<1+{\textstyle \frac{2}{7}}\epsilon\quad\mbox{for all }t.\label{eq:F*}\end{equation}

Now use lemma \ref{lem:nologo} with $f_{\textrm{lemma \ref{lem:nologo}}}=G^{*}$,
$g_{\textrm{lemma \ref{lem:nologo}}}=F^{*}$ and the set $\mathcal{B}^{*}:=\{ t:|G^{*}(t)|\leq\mu_{i}-\Osc(g)-\frac{3}{7}\epsilon\}$
and again $\epsilon_{\textrm{lemma }\ref{lem:nologo}}=\frac{1}{7}\epsilon$.
The resulting functions are our $G$ and $F$. As above, from the
definition of $c_{2}$ (\ref{eq:defc2}), (\ref{eq:reqc2}), (\ref{eq:gi-g})
and (\ref{eq:gG*}) we get that $1-c_{1}<|G^{*}|<1+c_{1}$ so if $\delta<\frac{1}{4}c_{1}$
we may apply the lemma. Further,\[
|\mathcal{B}^{*}|\stackrel{(\ref{eq:gG*})}{\leq}|\{ t:|g_{i}|\leq\mu_{i}-\Osc(g)-{\textstyle \frac{2}{7}}\epsilon\}|\stackrel{(\ref{eq:giosc})}{<}\delta\]
so that $\left\Vert F-F^{*}\right\Vert _{\infty}\leq C\sqrt{\delta}$
and if $\delta$ is sufficiently small, $\left\Vert F-F^{*}\right\Vert _{\infty}\leq\frac{1}{7}\epsilon$
so with (\ref{eq:F*}) we have what we need for $F$. For $G$ we
get \[
\mu_{i}-\Osc(g)-{\textstyle \frac{4}{7}}\epsilon<|G(t)|<\mu_{i}+{\textstyle \frac{3}{7}}\epsilon\]
so $\Osc(G)<\Osc(g)+\epsilon$. Fixing $\delta$ sufficiently small
and $i$ sufficiently large to satisfy all the past requirements the
lemma is done.
\end{proof}
\begin{lem}
\label{lem:killf}Let $f$, $g$ and $\epsilon$ be as in lemma \ref{lem:smooth}
but with the additional requirement \[
1-c_{3}<|f|<1\quad1-c_{3}<|g|<1+c_{3}.\]
 Then one can extend $f$ and $g$ to $F$ and $G$ such that
\begin{enumerate}
\item \label{enu:kill F}$1-\epsilon<|F(t)|<1+\epsilon$ for all $t$,
\item $\left\Vert F-f\right\Vert _{\infty}<C\sqrt{\left\Vert 1-|f|\,\right\Vert _{\infty}}$,
\item \label{enu:kill G}$\Osc(G)<\Osc(g)+\epsilon$,
\item $\left\Vert G-g\right\Vert _{\infty}<C\sqrt{\left\Vert 1-|f|\,\right\Vert _{\infty}}$.
\end{enumerate}
\end{lem}
\begin{proof}
This follows easily by applying lemma \ref{lem:halff} repeatedly
(to preserve the requirement $|f|<1$ you need to multiply $f$ and
$g$ by normalization factors). Let us do it in detail nonetheless.
Denote $\rho=\left\Vert 1-|f|\,\right\Vert _{\infty}$ and let $\delta$
be some parameter sufficiently small to be fixed later. We define
$f_{0}=f$ and $g_{0}=g$ and then inductively $f_{i}$ and $g_{i}$
satisfying the following properties
\begin{enumerate}
\item \label{enu:extends}$|\widehat{f_{i}}(n)|=|\widehat{g_{i}}(n)|$ for
all $n$. $(1+\delta2^{-i})f_{i}$ extends $f_{i-1}$ and $(1+\delta2^{-i})g_{i}$
extends $g_{i-1}$.
\item \label{enu:fi1rho1}$1-(\rho+2i\delta)2^{-i}<|f_{i}|<1$.
\item $\left\Vert f_{i}-f_{i-1}\right\Vert _{\infty}\leq C\sqrt{\rho2^{-i}}$
for $i>0$
\item \label{enu:muimui-1}$\Osc(g_{i})<\Osc(g)+\delta(1-2^{-i})$
\item $\left\Vert g_{i}-g_{i-1}\right\Vert _{\infty}\leq C\sqrt{\rho2^{-i}}$
for $i>0$
\end{enumerate}
Let us verify that the induction holds. At the $i$'th step ($i\geq1$)
we wish to apply lemma \ref{lem:halff} to the functions $f_{i-1}$
and $g_{i-1}$. For this to hold we must have that (\ref{eq:reqc2})
holds for our functions $f_{i-1}$ and $g_{i-1}$. First we note that
if $c_{3}<c_{2}$ and $\delta<\frac{1}{3}c_{2}$ then the second induction
assumption assures us that $1-c_{2}<|f_{i-1}|<1$. As for $g_{i-1}$,
we first estimate $\left\Vert g_{i-1}\right\Vert _{\infty}$ by noting
that \begin{align*}
\left\Vert g_{i-1}\right\Vert _{\infty} & \leq\left\Vert g_{i-1}\right\Vert _{2}+\Osc(g_{i-1})\stackrel{\textrm{\ref{enu:muimui-1}}}{\leq}\left\Vert g_{i-1}\right\Vert _{2}+\Osc(g)+\delta=\\
 & \stackrel{\textrm{\ref{enu:extends}}}{=}\left\Vert f_{i-1}\right\Vert _{2}+\Osc(g)+\delta\stackrel{\textrm{\ref{enu:fi1rho1}}}{\leq}1+2c_{3}+\delta\end{align*}
A similar calculation shows that\[
\left\Vert g_{i-1}\right\Vert _{\infty}\geq1-(\rho+2i\delta)2^{-i}-2c_{3}-\delta\geq1-3c_{3}-2\delta.\]
With these estimates we write \[
\left\Vert \,|g_{i-1}|-1\right\Vert _{\infty}\leq\Osc(g_{i-1})+\left|\,\left\Vert g_{i-1}\right\Vert _{\infty}-1\right|\leq5c_{3}+3\delta.\]
Hence if we define $c_{3}:=\frac{1}{10}c_{2}$ and ensure $\delta<\frac{1}{6}c_{2}$
then the requirements for lemma \ref{lem:halff} are assured.

We now apply lemma \ref{lem:halff} with $\epsilon_{\textrm{lemma }\ref{lem:halff}}=\delta2^{-i}$.
Call the resulting functions $f^{*}$ and $g^{*}$ and define $f_{i}=f^{*}/(1+\delta2^{-i})$
and $g_{i}=g^{*}/(1+\delta2^{-i})$. It is quite easy to verify that
all the inductive assumptions hold --- let us do two examples in detail.
First, let us verify the left hand side of \ref{enu:fi1rho1}. We
have\[
|f_{i}|\geq\frac{1-\frac{1}{2}\left\Vert 1-|f_{i-1}|\,\right\Vert -\delta2^{-i}}{1+\delta2^{-i}}\stackrel{\textrm{\ref{enu:fi1rho1}}}{>}\frac{1-(\rho+(2i-1)\delta)2^{-i}}{1+\delta2^{-i}}>1-(\rho+2i\delta)2^{-i}\]
where the inequality labeled \ref{enu:fi1rho1} uses this clause for
$i-1$ inductively. Secondly discuss \ref{enu:muimui-1}. We have
$\Osc(g_{i})=\Osc(g^{*})/(1+\delta2^{-i})<\Osc(g^{*})$ while clause
\ref{enu:G mu g nu} of lemma \ref{lem:halff} gives $\Osc(g^{*})<\Osc(g_{i-1})+\delta2^{-i}$.
This completes the induction.

Now take $i$ sufficiently large and define $\lambda_{i}=\prod_{j=1}^{i}(1+\delta2^{-j})$,
$F:=f_{i}\lambda_{i}$ and $G:=g_{i}\lambda_{i}$. We note that \[
|F|\leq\lambda_{i}\leq1+C\delta\qquad|F|\geq1-(\rho+2i)2^{-i}\geq1-C2^{-i}\]
hence for $i$ sufficiently large and $\delta$ sufficiently small
the first requirement from $F$ is satisfied. The other requirements
from $F$ and $G$ may be verified with similar ease.
\end{proof}

\begin{proof}
[Proof of the theorem]\label{page:start}Let $\varphi:\mathbb{T}\to\mathbb{T}$
be a function which does one rotation around $0$ very quickly, namely
for some $\epsilon>0$ to be fixed later define \[
\varphi(t)=\begin{cases}
e^{2\pi it/\epsilon} & t\leq\epsilon\\
1 & t>\epsilon.\end{cases}\]
Let $f_{1}$ be a trigonometric polynomial satisfying $\left\Vert f_{1}-\varphi\right\Vert _{\infty}<\epsilon$,
$\deg f_{1}<C\epsilon^{-C}$. Clearly $\wind f_{1}=1$ and $\left\Vert f_{1}-1\right\Vert _{2}<C\sqrt{\epsilon}$.
Examine the random function \[
h=\sum_{n\neq0}\pm\widehat{f_{1}(n)}e^{int}.\]
As before we have $\left\Vert h\right\Vert _{\infty}\leq C\sqrt{\epsilon\log1/\epsilon}$
with high probability. Pick one combination of signs $\xi_{n}$ such
that this holds and define $g_{1}=\widehat{f}_{1}(0)+\sum\xi_{n}\widehat{f}_{1}(n)e^{int}$.
Clearly $|\widehat{f}_{1}(n)|=|\widehat{g_{1}}(n)|$ and if $\epsilon$
is sufficiently small $\wind g_{1}=0$.

We now apply lemma \ref{lem:killf} inductively as follows. For the
even $j$-s we apply it with \[
f_{\textrm{lemma \ref{lem:killf}}}=f_{j-1}\nu_{j}\qquad g_{\textrm{lemma \ref{lem:killf}}}=g_{j-1}\nu_{j}\qquad\nu_{j}=\frac{1}{\max\{\left\Vert f_{j-1}\right\Vert _{\infty},1\}}\]
 and then define $f_{j}=F_{\textrm{lemma \ref{lem:killf}}}$ and $g_{j}=G_{\textrm{lemma \ref{lem:killf}}}$.
For the odd $j$-s we reverse the roles of $f$ and $g$ i.e.~take
the normalization factor to be $\nu_{j}=1/\max\{||g_{j-1}||_{\infty},1\}$
and then $f_{\textrm{lemma \ref{lem:killf}}}=\nu_{j}g_{j-1}$ and
$g_{\textrm{lemma \ref{lem:killf}}}=\nu_{j}f_{j-1}$ etc. In both
cases we take $\epsilon_{\textrm{lemma \ref{lem:killf}}}=2^{-j}\epsilon$.
An argument similar to that of lemma \ref{lem:killf} now shows that
throughout this process \begin{equation}
\Osc(f_{j})<5\epsilon2^{-j}\quad\Osc(g_{j})<5\epsilon2^{-j}.\label{eq:fjgjconst}\end{equation}
Let us show (\ref{eq:fjgjconst}) for the case of $j$ even (the other
case is identical). Here $f_{j}$ follows immediately since by clause
\ref{enu:kill F} of lemma \ref{lem:killf}, $\Osc(f_{j})<2\epsilon2^{-j}$.
As for $g$, the same clause \ref{enu:kill F} in step $j-1$ shows
that $\Osc(g_{j-1})<4\epsilon2^{-j}$ so $\Osc(g_{j-1}\nu_{j})<4\epsilon2^{-j}$.
We apply clause \ref{enu:kill G} of lemma \ref{lem:killf} in step
$j$ and get, $\Osc(g_{j})<5\epsilon2^{-j}$.

A similar argument shows that both are close to $1$ in the sense
that \begin{equation}
\left\Vert \,|f_{j}|-1\right\Vert _{\infty}<6\epsilon2^{-j}\quad\left\Vert \,|g_{j}|-1\right\Vert _{\infty}<6\epsilon2^{-j}.\label{eq:near1}\end{equation}
Again we demonstrate this under the assumption that $j$ is even.
For $f_{j}$ this is immediate from clause \ref{enu:kill F} of lemma
\ref{lem:killf}. For $g_{j}$, since $\left\Vert g_{j}\right\Vert _{2}=\left\Vert f_{j}\right\Vert _{2}$
we get that $\left|\,\left\Vert g_{j}\right\Vert _{2}-1\right|<\epsilon2^{-j}$
and since $\left|g_{j}(t)-\left\Vert g_{j}\right\Vert _{2}\,\right|\leq\Osc(g_{j})$
we get  (\ref{eq:near1}).

This implies that if $\epsilon$ is sufficiently small the induction
actually works in the sense that $1-c_{3}<\nu_{j}f_{j-1}<1+c_{3}$
and $1-c_{3}<\nu_{j}g_{j-1}<1+c_{3}$ are preserved throughout. Further
it implies$\left\Vert f_{j+1}-f_{j}\right\Vert _{\infty}<C\sqrt{\epsilon2^{-j}}$
and in particular we get that (if $\epsilon$ is sufficiently small)
that $\wind f_{j}=\wind f_{j+1}$, that $f=\lim f_{j}$ exists and
is continuous, that $|f(t)|=1$ for all $t$ and that $\wind f=\wind f_{1}=1$.
Similarly we get that $g=\lim g_{j}$ exists and is continuous, that
$|g(t)|=1$ for all $t$ and that $\wind g=0$. The property that
$|\widehat{f_{j}}(n)|=|\widehat{g_{j}}(n)|$ is preserved through
the limit so $|\widehat{f}(j)|=|\widehat{g}(j)|$ for all $j$ and
the theorem is proved. 
\end{proof}

\end{document}